\documentclass[12pt]{article}

\usepackage[dvips]{color}
\usepackage{amsfonts}
\usepackage{graphicx}
\usepackage{amsmath}
\usepackage{amssymb}
\usepackage{amsthm}
\usepackage{mathrsfs}
\usepackage{mathrsfs}
\newtheorem{theorem}{Theorem}[section]

\theoremstyle{definition}

\theoremstyle{remark}

\numberwithin{equation}{section}

\newfam\msbfam
\font\tenmsb=msbm10  scaled \magstep1 \textfont\msbfam=\tenmsb
\font\sevenmsb=msbm7 scaled \magstep1 \scriptfont\msbfam=\sevenmsb
\font\fivemsb=msbm5  scaled \magstep1 \scriptscriptfont\msbfam=\fivemsb
\def\Bbb{\fam\msbfam \tenmsb}

\def\CC{{\Bbb C}}

\def\ZZ{{\Bbb Z}}

\def\ss{\subseteq}
\def\ra{\rightarrow}

 \def\HollowBox #1#2{{\dimen0=#1 \advance\dimen0 by -#2       
       \dimen1=#1 \advance\dimen1 by #2                       
        \vrule height #1 depth #2 width #2                    
        \vrule height 0pt depth #2 width #1                   
        \llap{\vrule height #1 depth -\dimen0 width \dimen1}%
       \hskip -#2                                             
       \vrule height #1 depth #2 width #2}}                   
 \def\BoxOpTwo{\mathord{\HollowBox{6pt}{.4pt}}\;}             

\def\endpf{\hfill $\BoxOpTwo$ \smallskip \\ }

\def\Aut{\hbox{Aut}\,}
\def\dbar{\overline{\partial}}

\newfam\msbbfam
\font\tenmsbb=msbm10  scaled \magstep1 \textfont\msbbfam=\tenmsbb
\font\sevenmsbb=msbm7  scaled \magstep1 \scriptfont\msbbfam=\sevenmsbb
\font\fivemsbb=msbm5    scaled \magstep1 \scriptscriptfont\msbbfam=\fivemsbb

\begin{document}

\begin{center}
\Large \bf Semicontinuity of the Automorphism Groups of Domains with Rough Boundary\footnote{{\bf Subject 
Classification Numbers:}  32M05, 32M17, 32M25, 30F10.}
\end{center}
\vspace*{.12in}

\begin{center}
Steven G. Krantz
\end{center}

\date{\today}

\begin{quote}
{\bf Abstract:}  Based on some ideas of Greene and Krantz, we study the semicontinuity of automorphism
groups of domains in one and several complex variables.  We show that semicontinuity
fails for domains in $\CC^n$, $n > 1$, with Lipschitz boundary, but it holds for domains
in $\CC^1$ with Lipschitz boundary.  Using the same ideas, we develop some other
concepts related to mappings of Lipschitz domains.  These include Bergman curvature, stability properties
for the Bergman kernel, and also some ideas about equivariant embeddings.
\end{quote}
\vspace*{.25in}

\markboth{STEVEN G. KRANTZ}{SEMICONTINUITY OF AUTOMORPHISM GROUPS}

\section{Introduction}

A domain in $\CC^n$ is a connected open set.  If $\Omega$ is a domain,
then we let $\hbox{Aut}\, (\Omega)$ denote the group (under the binary
operation of composition of mappings) of biholomorphic self-maps of
$\Omega$.  When $\Omega$ is a bounded domain, $\hbox{Aut}\, (\Omega)$ is
a real (never a complex) Lie group.

A notable theorem of Greene/Krantz [GRK2] says the following:

\begin{theorem} \sl
Let $\Omega_0$ be a smoothly bounded, strongly pseudoconvex domain with defining
function $\rho_0$ (see [KRA1] for the concept of defining function).   There is an $\epsilon > 0$ so that, if $\rho$ is a defining
function for a smoothly bounded, strongly pseudoconvex domain $\Omega$ with
$\|\rho_0 - \rho\|_{C^k} < \epsilon$ (some large $k$) then the automorphism
group of $\Omega$ is a subgroup of the automorphism group of $\Omega_0$.
Furthermore, there is a diffeomorphism $\Phi: \Omega \ra \Omega_0$ such
that the mapping
$$
\hbox{Aut}\, (\Omega) \ni \varphi \longmapsto \Phi \circ \varphi \circ \Phi^{-1}
$$
is an injective group homomorphism of $\hbox{Aut}\, (\Omega)$ into $\hbox{Aut}\, (\Omega_0)$.
\end{theorem}

In what follows we shall refer to this result as the ``semicontinuity theorem.''

It should be noted that, although this theorem was originally proved for strongly
pseudoconvex domains in $\CC^n$, the very same proof shows that the
result is true in $\CC^1$ for any smoothly bounded domain $\Omega_0$.  In fact
the proof, while parallel to the original proof in [GRK2], is considerably
simpler in the one-dimensional context.

The original proof of this result, which was rather complicated, used stability results for the Bergman
kernel and metric established in [GRK1] and also the idea of Bergman
representative coordinates.   An alternative approach, using normal families,
was developed in [KIM].	  The paper [GRK3] produced a method for
deriving a semicontinuity theorem when the domain boundaries
are only $C^2$.  The more recent work [GKKS] gives a new
and more powerful approach to this matter of reduced boundary smoothness.
The paper [KRA2] gives yet another approach to the matter, and proves
a result for finite type domains.

It is geometrically natural to wonder whether there is a semicontinuity
theorem when the boundary has smoothness of degree less than 2.	  On the
one hand, experience in geometric analysis suggests that $C^2$ is a natural
cutoff for many positive results (see [KRP]).  On the other hand,
Lipschitz boundary is very natural from the point of view of dilation and
other geometric operations.

The purpose of this paper is to show that the semicontinuity theorem fails
for domains in $\CC^n$, $n > 1$, with Lipschitz boundary.  But it holds
for domains in $\CC^1$ with Lipschitz boundary.  The reason for this
difference is connected, at least implicitly, with the failure of the
Riemann mapping theorem in several complex variables.  We shall explain
this point in more detail as the presentation develops.

\section{The Several-Complex-Variable Situation}

The main result of this section is the following:

\begin{theorem} \sl
Let $n > 1$ and consider domains in $\CC^n$.
There is a sequence $\Omega_j$ of strongly pseudoconvex domains
with Lipschitz boundary and another domain  $\Omega$ with
Lipschitz boundary so that $\Omega_j \ra \Omega$ in the Lipschitz
topology on defining functions and so that
\begin{enumerate}
\item[{\bf (a)}]  For each $j$, $\Aut(\Omega_j) = \ZZ$;
\item[{\bf (b)}]  $\Aut(\Omega) = \{\hbox{id}\}$.
\end{enumerate}
\end{theorem}

See [HEL] for a consideration of strongly pseudoconvex domains
with less than $C^2$ boundary. This result shows that the
semicontinuity theorem fails for domains with Lipschitz
boundary.

It should be understood that all the domains considered in this 
paper have finite connectivity.  In particular, the complement of
the domain only has finitely many components.  And each component
of the complement has Lipschitz boundary.  We do {\it not} allow boundary
components that are a single point.  Each boundary component is the closure
of an open set.

We shall use some ideas in [LER] in constructing the
example enunciated in the theorem.   We shall make our
construction in $\CC^2$.  But it is easy to produce analogous
examples in any $\CC^n$.
\smallskip \\

\noindent {\bf Proof of the Theorem:}	Let $\psi \in C_c^\infty(\CC^n)$
be such that
\begin{enumerate}
\item[{\bf (i)}]  $\hbox{supp}\, \psi \ss B(0,1)$;
\item[{\bf (ii)}]  $\psi \geq 0$;
\item[{\bf (iii)}]  $\psi(0) = 1$.
\end{enumerate}

We will build our domains by modifying the unit ball $B$ in $\CC^2$.  We will
make particular use of these automorphisms of the unit ball, for $a$ a complex
number of modulus less than 1:
$$
\Psi_a(z_1, z_2) = \left ( \frac{z_1 - a}{1 - \overline{a} z_1} , \frac{\sqrt{1 - |a|^2} z_2}{1 - \overline{a} z_1} \right ) \, .
$$
See [RUD].

We define 
$$
\eta_1(z_1, z_2) = - 1 + |z_1|^2 + |z_2|^2 - (1/10) \psi \biggl ( 10 \bigl ( (z_1, z_2) - (\sqrt{3/4},1/2) \bigr ) \biggr ) \, .
$$
Set
$$
U_1 = \{(z_1, z_2) \in \CC^2: \eta_1(z_1, z_2) < 0\} \, .
$$
Clearly $U_1$ is a domain with smooth boundary.  It is a ball with
a ``bump'' attached at the point $(\sqrt{3/4},1/2)$.  

Now define
$$
\Omega_1 = \bigcup_{j= -\infty}^\infty \Psi^{2^j}_{1/10} (U_1) \, .
$$
We see that $\Omega_1$ has infinitely many bumps which accumulate at
the points $(1,0)$ and $(-1,0)$.  It is because of those accumulation
points that the boundary of $\Omega_1$ is only Lipschitz.
						
In general we let, for $k \geq 2$,  
$$
\eta_k(z_1, z_2) = - 1 + |z_1|^2 + |z_2|^2 - (1/10^k) \psi \biggl ( 10^k \bigl ( (z_1, z_2) - (\sqrt{(1/2)^{k-1} - (1/2)^{2k}} ,1- (1/2)^k) \bigr ) \biggr ) \, .
$$
Set
$$
U_k = \{(z_1, z_2) \in \CC^2: \eta_k(z_1, z_2) < 0\} \, .
$$
Clearly $U_k$ is a domain with smooth boundary.  It is a ball with
a ``bump'' attached at the point $(\sqrt{(1/2)^{k-1} - (1/2)^{2k} },1- (1/2)^k)$.  

Now define, for $k \geq 2$, 
$$
\Omega_k = \Omega_{k-1} \cup \bigcup_{j= -\infty}^\infty \Psi^{2^{j + k - 1}}_{1/10} (U_k) \, .
$$
We see that $\Omega_k$ has infinitely many bumps which accumulate at
the points $(1,0)$ and $(-1,0)$.  It is because of those accumulation
points that the boundary of $\Omega_1$ is only Lipschitz.

Finally we let
$$
\Omega = \bigcup_{k=1}^\infty \Omega_k \, .
$$

Now it is clear that $\Omega_k \ra \Omega$ in the Lipschitz topology
on defining functions.  Furthermore, the ideas in [LER] show that
the automorphism group of $\Omega_k$ consists precisely
of the mappings $\Psi_{1/10}^{2^{j + k - 1}}$, $j \in \ZZ$.  So the
automorphism group of $\Omega_k$ is canonically isomorphic
to $\ZZ$.  But it is also clear that the automorphism group
of $\Omega$ consists of the identity alone.

That completes the construction described in the theorem.
\endpf

\section{The One-Variable Situation}

The one-variable result is this:

\begin{theorem} \sl
Consider domains in $\CC^1$.
Let $\Omega_0 \ss \CC^1$ be a bounded domain with Lipschitz boundary
and defining function $\rho_0$.
If $\epsilon > 0$ is sufficiently small then, whenever $\Omega$
is a bounded domain with Lipschitz boundary and defining
function $\rho$ satisfying $\|\rho_0 - \rho\|_{\rm Lip} < \epsilon$
then the automorphism group of $\Omega$ is a subgroup of the automorphism
group of $\Omega_0$.  Moreover, there is a diffeomorphism $\Phi: \Omega \ra \Omega_0$
so that the mapping
$$
\Aut(\Omega) \ni \varphi \longmapsto \Phi \circ \varphi \circ \Phi^{-1} \in \Aut(\Omega_0)
$$
is an injective group homomorphism.
\end{theorem}

We see here that the situation is in marked contrast to that for several complex variables.
Our proof of this result will rely on uniformization for planar domains, a result
which has no analogue in several complex variables.  
\smallskip \\

\noindent {\bf Proof of the Theorem:}  Fix the domain $\Omega_0$ and let
$\Omega$ be of distance $\epsilon$ from $\Omega_0$ in the Lipschitz topology.

It is a standard result of classical function theory that a finitely
connected domain in the plane, with no component of the complement equal
to a point, is conformally equivalent to the plane with finitely many
nontrivial closed discs excised---see [AHL] or [KRA3]. Call this conformal mapping the
``normalization'' of the domain. What is particularly nice about this
result is that the proof is constructive and it is straightforward to see
that the normalization of $\Omega$ is close to the normalization of $\Omega_0$
just because $\Omega$ is close to $\Omega_0$.   Indeed the normalization of
$\Omega$ will be close to that of $\Omega_0$ in the $C^2$ topology.   Just because
once it is close in the Lipschitz topology then it is automatically close
in a smoother topology (because the boundary consists of finitely many nontrivial circles).

Thus we may apply the one-dimensional version of the semicontintuity theorem for $C^2$ boundary
to see that the automorphism group of the normalization of $\Omega$ is a subgroup
of the automorphism group of $\Omega_0$.  And the diffeomorphism $\Phi$ exists
as usual.   Now we may use the normalizing conformal mapping to transfer
this result back to the original domains $\Omega_0$ and $\Omega$.

That completes the proof.
\endpf

We note that another approach to constructing the normalization map is by way
of Green's functions.  This method is also quite explicit and constructive.
Stability results for elliptic boundary value problems are well known.
So this again leads to a proof of the semicontinuity theorem by transference
to the normalized domain.

\section{Related Results in One Complex Dimension}

Key to the work of Greene-Krantz in [GRK1] and [GKR2] is a stability
result for the Bergman kernel.  In that theorem, the authors
consider a base domain $\Omega_0$ and a ``nearby'' domain $\Omega$.
As usual, we define ``nearby'' in terms of closeness of the defining
functions in a suitable topology.   But it is useful to note that,
in this circumstance, there is a diffeomorphism $\Pi: \Omega \ra \Omega_0$ 
which is close to the identity in a suitable $C^k$ topology.  With
this thought in mind, Greene and Krantz proved the following:

\begin{theorem} \sl
Let $\Omega_0$ be a fixed, smoothly bounded, strongly pseudoconvex domain.
Let $\Omega$ be a domain which is ``$\epsilon$-close'' to $\Omega_0$ in a
$C^k$ topology.  Let $\Pi$ be the mapping described in the preceding paragraph.
If $\epsilon$ is small enough, then the Bergman kernel $K_\Omega$ for
$\Omega$ is close to $K_{\Omega_0} \circ \Pi$ in the $C^m$ topology
for some $0 < m < k$.
\end{theorem}

This result also holds in one complex dimension, and the proof in that context
is actually much easier.

Our remark now is that this theorem is actually true in the Lipschitz
topology.  We use the argument of the last section.  Namely, if $\Omega$
is close to $\Omega_0$ in the Lipschitz topology, then the normalization
of $\Omega$ is close to the normalization of $\Omega_0$ in a smooth
topology.  This the one-dimensional version of Theorem 4.1 applies
to the normalized domains.  The result follows.

\section{Equivariant Embeddings}

A lovely result of Maskit [MAS] is the following:

\begin{theorem} \sl
Let $\Omega \ss \CC$ be {\it any} planar domain.
Then there is a univalent, holomorphic embedding
$\Phi: \Omega \rightarrow \CC$ so that the automorphism
group of the image domain $\Phi(\Omega)$ consists only
of linear fractional transformations.
\end{theorem}

An elegant corollary of Maskit's result is that 
if $\varphi$ is any automorphism of a planar domain
that fixes three points then $\varphi$ is the identity
mapping.  This follows because it is clear that any linear
fractional transformation that fixes three points is
the identity.

We would like to remark here that the ideas in this paper
give a ``poor man's version'' of this theorem.  For
let $\Omega$ be any domain with Lipschitz boundary as we have
been discussing.  So each component of the complement is the
closure of a region having Lipschitz boundary.  Now the normalizing
map sends this domain $\Omega$ to a planar domain bounded by finitely many
disjoint circles.  It is easy to see, using Schwarz reflection and
Schwarz's lemma, that any conformal self-map of such a domain must
be linear fractional.   So any such map that fixes three points
must be the identity.

\section{The Bun Wong/Rosay Theorem}

A classical result in several complex variables is this (see [ROS], [WON]:

\begin{theorem} \sl
Let $\Omega \ss \CC^n$ be a bounded domain.  Let $P \in \partial \Omega$
and assume that $\partial \Omega$ is strongly pseudconvex in a neighborhood
of $P$.  Suppose that there are a point $X \in \Omega$ and automorphisms
$\varphi_j$ of $\Omega$ such that $\varphi_j(X) \ra P$ as $j \ra \infty$.
Then $\Omega$ is biholomorphic to the unit ball.
\end{theorem}

In a similar spirit, Krantz [KRA4] proved the following result:

\begin{theorem} \sl
Let $\Omega \ss \CC$ be a bounded domain and let $P \in \partial \Omega$
have the property that $\partial \Omega$ near $P$ is a $C^1$ curve.
Suppose that there are a point $X \in \Omega$ and automorphisms
$\varphi_j$ of $\Omega$ such that $\varphi_j(X) \ra P$ as $j \ra \infty$.
Then $\Omega$ is conformally equivalent to the unit disc.
\end{theorem}

In this section we will re-examine Theorem 6.2 in the context of this paper, that
is in relation to finitely connected domains with Lipschitz boundary.  As noted,
such a domain is conformally equivalent to a domain $\widehat{\Omega}$ whose boundary consists
of finitely many circles.  Now we have the following possibilities:
\begin{enumerate}
\item[{\bf (a)}]  If $\partial \widehat{\Omega}$ consists of just one circle, then
$\widehat{\Omega}$ is the disc, and there is nothing to prove.
\item[{\bf (b)}]  If $\partial \widehat{\Omega}$ consists of two circles, one
inside the other, then $\widehat{\Omega}$ is (conformally equivalent to) an annulus.
Then the automorphism group of such a domain is two copies of the unit circle.
In particular, it is compact.  So the hypotheses of Theorem 6.2 cannot obtain.
\item[{\bf (c)}]  If $\partial \widehat{\Omega}$ consists of two circls, neither
of which is inside the other, then the domain is unbounded.  The automorphism
group of such a domain is compact, and the hypotheses of Theorem 6.2 do not apply.
\item[{\bf (d)}]   If $\partial \widehat{\Omega}$ consists of at least three circles,
with all the circles but one lying inside the other one, then it is well known (see [JUL] or [HEI2])
that the automorphism group of $\widehat{\Omega}$ is finite.  Then the hypotheses
of Theorem 6.2 cannot obtain.
\end{enumerate}

Thus we see by inspection that Theorem 6.2 is true in the context of the domains
that we have been discussing in this paper. 

\section{Curvature of the Bergman Metric}

It is a matter of considerable interest to know the curvature properties of
the Bergman metric on a planar domain.  In particular, negativity of the curvature
near the boundary is a useful analytic tool (see [GRK1]).  If $\Omega$
is a planar domain with Lipschitz boundary, then its normalized domain
is bounded by finitely many circles.  The asymptotic boundary behavior
of the Bergman kernel on such a domain is very well understood---see [APF].
In particular, the kernel near a boundary point $P$ is asymptotically very much like
the kernel for the disc.   Thus a straightforward calculation confirms
that the curvature of the Bergman metric near the boundary is negative.
Of course this statement pulls back to the original domain in a natural
way.

\section{Closing Remarks}

It is natural to want to consider the results presented here in either
the $C^1$ topology or even the $C^{2-\epsilon}$ topology.  At this
time the techniques are not available to attack those questions.

In several complex variables, one would also like to prove semicontinuity
theorems for broad classes of domains.  This will be the subject for future
papers.

\newpage

\noindent {\Large \sc References}
\vspace*{.2in}

\begin{enumerate}

\item[{\bf [APF]}]  L. Apfel, Localization Properties and Boundary
Behavior of the Bergman Kernel, thesis, Washington University
in St.\ Louis, 2003.

\item[{\bf [GKKS]}] R. E. Greene, K.-T. Kim, S. G. Krantz, and
A.-R. Seo, Semi-continuity of automorphism groups of strongly
pseudoconvex domains: the low differentiability case,
preprint.

\item[{\bf [GRK1]}] R. E. Greene and S. G. Krantz, Deformations
of complex structure, estimates for the $\dbar$-equation, and
stability of the Bergman kernel, {\it Advances in Math.}
43(1982), 1--86.

\item[{\bf [GRK2]}] R. E. Greene and S. G. Krantz, The
automorphism groups of strongly pseudoconvex domains, {\it
Math. Annalen} 261(1982), 425-446.

\item[{\bf [GRK3]}] R. E. Greene and S. G. Krantz, Normal
families and the semicontinuity of isometry and automorphism
groups, {Math.\ Z.} 190(1985), 455--467.

\item[{\bf [HEI1]}]  M. Heins, A note on a theorem of Rad\'{o}
concering the $(1,m)$ conformal maps of a multiply connected
region into itself, {\it Bulletin of the AMS} 47(1941),
128--130.

\item[{\bf [HEI2]}] M. Heins, On the number of 1-1 directly
conformal maps which a multiply-connected plane regions of
finite connectivity $p (> 2)$ admits onto itself, {\it
Bulletin of the AMS} 52(1946), 454--457.

\item[{\bf [HEL]}] G. M. Henkin and J. Leiterer, {\it Theory of
Functions on Strictly Pseudoconvex Sets with Nonsmooth
Boundary}, with German and Russian summaries, Report MATH
1981, 2. Akademie der Wissenschaften der DDR, Institut f\"{u}r
Mathematik, Berlin, 1981.

\item[{\bf [HES]}]  Z.-X. He and O. Schramm, Fixed points, Koebe uniformization and circle
packings, {\it Ann.\ of Math.} 137(1993), 369--406.

\item[{\bf [JUL]}] G. Julia, Le\c{c}ons sur la repr\'{e}sentation
conforme des aires multiplement connexes, Paris, 1934.

\item[{\bf [KIM]}] Y. W. Kim, Semicontinuity of compact
group actions on compact differentiable manifolds, {\it Arch.\
Math.} 49(1987), 450--455.

\item[{\bf [KRA1]}]  S. G. Krantz, {\it Function Theory of Several
Complex Variables}, 2nd ed., American Mathematical Society,
Providenc, RI, 2001.

\item[{\bf [KRA2]}] S. G. Krantz, Convergence of automorphisms
and semicontinuity of automorphism groups, {\it Real Analysis
Exchange}, to appear.

\item[{\bf [KRA3]}] S. G. Krantz, {\it Cornerstones of Geometric Function
Theory: Explorations in Complex Analysis}, Birkh\"{a}user Publishing,
Boston, 2006.

\item[{\bf [KRA4]}] S. G. Krantz, Characterizations of smooth
domains in $\CC$ by their biholomorphic self maps, {\it Am.\
Math.\ Monthly} 90(1983), 555--557.

\item[{\bf [KRP]}]  S. G. Krantz and H. R. Parks, {\it The Geometry
of Domains in Space}, Birkh\"{a}user Publishing, Boston, 1996.

\item[{\bf [LER]}] L. Lempert and L. Rubel, An independence
result in several complex variables, {\it Proc.\ Amer.\ Math.\
Soc.} 113(1991), 1055--1065.

\item[{\bf [MAS]}] B. Maskit, The conformal group of a plane
domain,  {\it Amer.\ J.\ Math.}, 90 (1968), 718--722.

\item[{\bf [ROS]}] J.-P. Rosay, Sur une characterization de la
boule parmi les domains de $\CC^n$ par son groupe
d'automorphismes, {\it Ann. Inst. Four. Grenoble} XXIX(1979),
91--97.

\item[{\bf [RUD]}]  W. Rudin, {\it Function Theory in the Unit
Ball of $\CC^n$}, Springer-Verlag, New York, 1980.

\item[{\bf [WON]}] B. Wong, Characterizations of the ball in
$\CC^n$ by its automorphism group, {\it Invent.\ Math.}
41(1977), 253--257.

\end{enumerate}
\vspace*{.25in}

\begin{quote}
Steven G. Krantz  \\
Department of Mathematics \\
Washington University in St. Louis \\
St.\ Louis, Missouri 63130  \\
{\tt sk@math.wustl.edu}
\end{quote}

\end{document}